\documentclass{amsart}

\usepackage{graphics}
\usepackage{epsfig}
\usepackage{amsmath}
\usepackage{amssymb,fancyhdr,txfonts,pxfonts}

\usepackage{paralist}

\usepackage{pstricks,pst-text,pst-grad,pst-node,pst-3dplot,pstricks-add,pst-poly,pst-coil} 
\usepackage{pst-fun,pst-blur} 

\usepackage{hyperref}

\newtheorem{theorem}{Theorem}
\theoremstyle{plain}

\newtheorem{corollary}{Corollary}

\newtheorem{definition}{Definition}

\newtheorem{lemma}{Lemma}

\newtheorem{proposition}{Proposition}

\begin{document}

\title[Skew-Field of Trace-Preserving Endomorphisms in Affine Plane]{Skew-Field of Trace-Preserving Endomorphisms, of Translation Group in Affine Plane}

\author[Orgest ZAKA]{Orgest ZAKA}
\address{Orgest ZAKA: Department of Mathematics-Informatics, Faculty of Economy and Agribusiness, 
Agricultural University of Tirana, Tirana, Albania.}
\email{ozaka@ubt.edu.al, gertizaka@yahoo.com, ozaka@risat.org}



\begin{abstract}
In this paper we will show how to constructed an Skew-Field with trace-preserving endomorphisms of the affine plane. Earlier in my paper, we doing a detailed description of endomorphisms algebra and trace-preserving endomorphisms algebra in an affine plane, and we have constructed an associative unitary ring for which trace-preserving endomorphisms. In this paper we formulate and prove an important Lemma, which enables us to construct a particular trace-preserving endomorphism, with the help of which we can construct the inverse trace-preserving endomorphisms of every trace-preserving endomorphism. At the end of this paper we have proven that the set of trace-preserving endomorphisms together with the actions of 'addition' and 'composition' (which is in the role of 'multiplication') forms a skew-field.

\end{abstract}


\maketitle

\section{Introduction}

The foundations for the study of the connections between axiomatic geometry and algebraic structures were set forth by D. Hilbert \cite{Hilbert1959geometry}, recently elaborated and extended in terms of the algebra of affine planes in, for example, \cite{Kryftis2015thesis}, [3, §IX.3, p.574]. Also great contributions in this direction have been made by, E.Artin in \cite{11},  H. S. M. Coxeter, in \cite{14},  Marcel Berger in \cite{Marcel.Berger.Geometry1.and.Geometry2}, Robin Hartshrone in \cite{Hartshorne.Foundations}, etc. Even earlier, in my works '\cite{1}, \cite{2}, \cite{3}, \cite{4}, \cite{5}, \cite{6}, \cite{7}, \cite{8}, \cite{15}, \cite{Zaka_2016}' I have brought up quite a few interesting facts about the association of algebraic structures with affine planes and with 'Desargues affine planes', and vice versa. 

In paper \cite{5}, we have done a detailed, careful description, of translations and dilation’s in affine planes. We have proven that the set of dilation’s regarding the composition action is a group, and set of translations is a commutative group. We have proved that translation group is a normal subgroup of the group of dilation’s. 

Even earlier in the paper \cite{3}, \cite{4}, \cite{5}, \cite{6}, we have shown ways how to build an algebraic structure on an affine plane and vice versa. We have shown in the paper \cite{3}, \cite{4}, how a skew-field is built over any line of an 'Descriptive affine plan', apparently bringing geometric illustrations to each construction. 

In this paper we do the reverse, apparently studying the endomorphisms of translation groups in affine planes, separating them from the trace-preserving endomorphisms, for which in the paper \cite{8}, we doing a detailed description of endomorphisms algebra and trace-preserving endomorphisms algebra in an affine plane, and we have constructed an associative unitary ring for which trace-preserving endomorphisms, together with two actions in the' addition 'and' composition '. Continuing our research in this direction, in this paper we will show how to constructed an skew-Field with trace-preserving endomorphisms of translation groups in the affine plane. In this paper we formulate and prove an important Lemma, which enables us to construct a particular trace-preserving endomorphism for translations group in affine plane, with the help of which we can construct the inverse trace-preserving endomorphisms of every trace-preserving endomorphism of translations group in affine plane.

\section{Preliminaries}
Let $\mathcal{P}$ be a nonempty set, which is called set of points, $\mathcal{L}$ a nonempty subset of $\mathcal{P}$, which is called set of lines, and an incidence relation $\mathcal{I}\subseteq \mathcal{P} \times \mathcal{L}$.  For a point $P \in \mathcal{P}$ and an line $\ell \in \mathcal{L}$, the fact $(P,\ell )\in \mathcal{I}$, (equivalent to $P\mathcal{I}\ell $) we mark $P\in \ell$ and read \textit{point} $P$ \textit{ is incident with a line} $\ell$ or a line passes through points $P$ (contains point $P$).

\begin{definition} \label{def.affine.plane} 
\cite{11},\cite{13},\cite{3} Affine plane is called the incidence
structure  $\mathcal{A}=(\mathcal{P},\mathcal{L},\mathcal{I})$ that satisfies
the following axioms:
\begin{description}
	\item[A.1] For every two different points $P$ and $Q$ $\in $ $\mathcal{P}$, there
exists exactly one line $\ell \in $ $\mathcal{L}$ incident with that points.
  \item[A.2] For a point $P$ $\in $ $\mathcal{P}$, and an line $\ell \in $ $\mathcal{L}$
such that $(P,\ell )\notin $ $\mathcal{I}$, there exists one and only one
line $r\in $ $\mathcal{L}$, incident with the point $P$ and such that $\ell
\cap r=\varnothing .$
  \item[A.2] In $\mathcal{A}$ there are three non-incident points with a line.
\end{description}
\end{definition}

Whereas a line of the affine plane we consider as sets of points of affine
plane with her incidents. Axiom A.1 implicates that tow different lines of $
\mathcal{L}$ many have a common point, in other words tow different lines of 
$\mathcal{L}$ either have no common point or have only one common point.

\begin{definition}
Two lines $ \ell , m \in \mathcal{L}$  that are matching or do not have
common point are called \textbf{parallel} and in this case is written  $\ell
\parallel m,$ and when they have only one common point we say that they are
expected.
\end{definition}

\begin{definition} 
Let it be $\mathcal{A}=(\mathcal{P},\mathcal{L},\mathcal{I})$ an
affine plane and $\mathcal{S}$=$\{\psi :\mathcal{P}\rightarrow \mathcal{P}%
| $ where $\psi -$is bijection$\}$  set of bijections to set points $\mathcal{P}$ on yourself. Collineation of affine plane $\mathcal{A}$, called a bijection $\psi \in \mathcal{S}$, such that%
\begin{equation}
\forall \ell \in \mathcal{L},\psi \left( \ell \right) \in \mathcal{L},
\label{eq}
\end{equation}
\end{definition}

Otherwise, a collineation of the affine plane $\mathcal{A}$ is a bijection
of set $\mathcal{P}$ on yourself, that preserves lines.  It is known
that the set of bijections to a set over itself is a group on associated
with the binary action '$\circ$' of composition in it, which is known as
total group or symmetric groups.

\begin{definition} \label{def.fixpoint} \cite{5}
An point $P$ of the affine plane $\mathcal{A}$ called fixed point his
associated with a collineation $\delta ,$ if coincides with the image itself
 $\delta (P),$ briefly when, $P=\delta (P).$
\end{definition}

\begin{definition} \label{def.dilation} \cite{11}, \cite{5}
A \textbf{Dilation} of an affine plane $\mathcal{A}=(\mathcal{P},%
\mathcal{L},\mathcal{I})$ called a its collineation $\delta $ such that%
\begin{equation}
\forall P\neq Q\in \mathcal{P},\delta \left( PQ\right) \Vert PQ
\end{equation}
\end{definition}

Let it be $\mathbf{Dil}_{\mathcal{A}}=\left\{ \delta \in \mathbf{Col}_{%
\mathcal{A}}|\delta -\text{is a dilation of }\mathcal{A}\right\} $ the
dilation set of affine plane $\mathcal{A}=(\mathcal{P},\mathcal{L},%
\mathcal{I})$.

\begin{theorem} \cite{5}, \cite{15}
The dilation set $\mathbf{Dil}_{\mathcal{A}}$ \ of affine plane $\mathcal{A%
}$ forms a group with respect to composition $\circ .$
\end{theorem}
\proof See, Theorem 2.4, in \cite{5}. \qed

\begin{definition} \label{def.trace} 
Let it be $\delta $ an dilation of affine plane $\mathcal{A}=(\mathcal{P},%
\mathcal{L},\mathcal{I})$, and $P$ a point in it. Lines that passes by points $P$ and $%
\delta (P),$ called \textbf{trace} of points $P$ regarding dilation $\delta .$
\end{definition}

Every point of a traces of a not-fixed point, to an affine plane
associated with its dilation has its own image associated with that dilation
in the same traces (see \cite{5}). We also know the result: If an affine plane $\mathcal{A}=(\mathcal{P},\mathcal{L},\mathcal{I}),$ has
two fixed points about an dilation then he dilation is \textit{identical
dilation} $id_{\mathcal{P}}$ of his (Theorem 2.12 in \cite{5}).

In an affine plane related to dilation $\delta \neq id_{\mathcal{P}}$ all
traces \ $P\delta \left( P\right) $ for all $P\in \mathcal{P}$ , or cross
the by a single point, or are parallel between themselves (see \cite{15}, \cite{5}).

\begin{definition} \cite{11}, \cite{5}
The \textbf{Translation} of an affine plane $\mathcal{A}=(%
\mathcal{P},\mathcal{L},\mathcal{I}),$ called the his identical dilation $id_{
\mathcal{P}}$ and every other of it's dilation, about which that the affine
plane has not fixed points.
\end{definition}

If $\sigma $ is an translation different from identical translation $id_{%
\mathcal{P}},$ then, all traces related to $\sigma $ form the
a set of parallel lines.

\begin{definition} \label{dif.direction}
For one translation $\sigma \neq id_{\mathcal{P}},$  the parallel equivalence classes of
the cleavage $\pi =\mathcal{L}/\Vert ,$ which contained tracks by  $\sigma $
of points of the plane $\mathcal{A}=(\mathcal{P},\mathcal{L},\mathcal{I})$
called the direction of his translation $\sigma $\ and marked with $\pi
_{\sigma }.$
\end{definition}

So, for $\sigma \neq id_{\mathcal{P}},$ the direction $\pi _{\sigma
}$ represented by single the trace (which is otherwise called, representative of direction) by  $\sigma $ for every point $P\in 
\mathcal{P},$ for translation $id_{\mathcal{P}}$  has
undefined direction.

Let it be $\alpha :\mathbf{Tr}_{\mathcal{A}}\longrightarrow \mathbf{Tr}_{%
\mathcal{A}},$ an whatever map of $\mathbf{Tr}_{\mathcal{A}},$ \ on
yourself. For every translation $\sigma,$ its image $\alpha \left(\sigma
\right)$ is again an translation, that can be $\alpha \left( \sigma \right)
=id_{\mathcal{P}}$ or $\alpha \left( \sigma \right) \neq id_{\mathcal{P}}$ .
So there is a certain direction or indefinite. The first equation, in
the case where \ $\sigma =id_{\mathcal{P}},$ takes the view $\alpha \left(
id_{\mathcal{P}}\right) =id_{\mathcal{P}},$ and the second $\alpha \left(
\sigma \right) \neq id_{\mathcal{P}},$ that it is not possible to $\alpha $
is map. To avoid this, yet accept that for every map $
\alpha :\mathbf{Tr}_{\mathcal{A}}\longrightarrow \mathbf{Tr}_{\mathcal{A}},$
is true this,
\[
\alpha \left( id_{\mathcal{P}}\right) =id_{\mathcal{P}}.
\]

\begin{proposition} \cite{5}
If translations $\sigma _{1}$ and $\sigma _{2}$ have the same direction with
translation $\sigma $ to an affine plane $\mathcal{A}=(\mathcal{P}
,\mathcal{L},\mathcal{I})$, then and composition $%
\sigma _{2}\circ \sigma _{1}$ has the same the direction, otherwise%
\begin{equation} \label{13}
\forall \sigma _{1},\sigma _{2},\sigma \in \mathbf{Tr}_{\mathcal{A}},\pi
_{\sigma _{1}}=\pi _{\sigma _{2}}=\pi _{\sigma }\Longrightarrow \pi _{\sigma
_{2}\circ \sigma _{1}}=\pi _{\sigma }.
\end{equation}
\end{proposition}

\begin{theorem} \cite{5}
Set $\mathbf{Tr}_{\mathcal{A}}$ of translations to an affine plane $\mathcal{%
A}$ form a group about the composition $\circ ,$ which is a sub-group of the
group $\left( \mathbf{Dil}_{\mathcal{A}},\circ \right) $ to dilations of
affine plane $\mathcal{A}$.
\end{theorem}

\begin{theorem} \cite{5} \label{thm.normal}
Group \ $\left( \mathbf{Tr}_{\mathcal{A}},\circ \right) $ of translations to
the affine plane $\mathcal{A}$ is normal sub- group of the group of
dilation’s \ $\left( \mathbf{Dil}_{\mathcal{A}},\circ \right) $ of him
plane.
\[
\forall \delta \in Dil_{\mathcal{A}},\forall \sigma \in Tr_{\mathcal{A}}\Rightarrow {\delta}^{-1} \circ \sigma \circ \delta \in Tr_{\mathcal{A}}.
\]
\end{theorem}

\begin{corollary} \cite{5} \label{cor.1}
For every dilation $\delta \in \mathbf{Dil}_{\mathcal{A}}$ and for every
translations $\sigma \in \mathbf{Tr}_{\mathcal{A}}$ of affine plane $%
\mathcal{A}=(\mathcal{P},\mathcal{L},\mathcal{I})
$, translations $\sigma $ and $\delta ^{-1}\circ \sigma \circ \delta $ of
his have the same direction.
\end{corollary}

\begin{corollary} \label{14}
The translations group $\left( \mathbf{Tr}_{\mathcal{A}},\circ \right) $
of an affine plane $\mathcal{A}$ is (Abelian) commutative Group.
\end{corollary}

In paper \cite{8}, we have defined the addition action for two maps $\alpha, \beta \in {(Tr_{\mathcal{A}})}^{Tr_{\mathcal{A}}}$, in set of maps, of commutative Group $(Tr_{\mathcal{A}}, \circ) $ of affine plane $\mathcal{A}$ in itself (see \cite{5}), so ${(Tr_{\mathcal{A}})}^{Tr_{\mathcal{A}}}=\left\{\alpha | \alpha: Tr_{\mathcal{A}}\rightarrow Tr_{\mathcal{A}}\right\}$, as below,

Let be $\alpha, \beta$ two different maps, such. Then for every 
\[
\forall \sigma \in Tr_{\mathcal{A}}  \Rightarrow \alpha(\sigma) \in Tr_{\mathcal{A}}, \text{ \ and \ } \beta(\sigma) \in Tr_{\mathcal{A}}, \] 
also,

\[ \forall \sigma \in Tr_{\mathcal{A}}  \Rightarrow [\alpha \circ \beta](\sigma)=\alpha (\beta(\sigma))\in Tr_{\mathcal{A}}. \] 

From the latter it turns out that the action of composition '$\circ$' in the set ${(Tr_{\mathcal{A}})}^{Tr_{\mathcal{A}}}$, is action induced by the action of composition '$\circ$' in the set $({\mathcal{P}})^{\mathcal{P}}$. If the associate translations $\sigma$, the unique translation  $\alpha (\sigma) \circ \beta(\sigma)$, obtained a new map $Tr_{\mathcal{A}}\rightarrow Tr_{\mathcal{A}}$, which we call \textbf{addition} of $\alpha$ with $\beta$, and mark with $\alpha + \beta$. 

\begin{definition} {\label{Def.Sum}} 
For every two maps  $\alpha, \beta \in {(Tr_{\mathcal{A}})}^{Tr_{\mathcal{A}}}$, the addition of them, that is marked $\alpha + \beta$, is called the map '${\alpha + \beta}:Tr_{\mathcal{A}} \rightarrow Tr_{\mathcal{A}}$', defined by,
\begin{equation} \label{15}
(\alpha + \beta)(\sigma)=(\alpha)(\sigma)\circ (\beta)(\sigma), \forall \sigma \in Tr_{\mathcal{A}}.
\end{equation}
\end{definition}

Accompanying any two maps $\alpha, \beta$ their sum $\alpha + \beta$, we obtain a new binary action in $Tr_{\mathcal{A}}$, that we call the addition of maps, of translations in the affine plane $\mathcal{A}$.

Thus obtained, algebra with two binary operations $({(Tr_{\mathcal{A}})}^{Tr_{\mathcal{A}}},+, \circ)$, where the sum of the two elements, whatsoever, its $\alpha, \beta$ in ${(Tr_{\mathcal{A}})}^{Tr_{\mathcal{A}}}$, is given by Definition \ref{Def.Sum}, and their composition is given by,

\begin{equation} \label{16}
\forall \sigma \in Tr_{\mathcal{A}}, [\alpha \circ \beta](\sigma)=\alpha (\beta(\sigma)).
\end{equation}

The Algebra $\left({(Tr_{\mathcal{A}})}^{Tr_{\mathcal{A}}},+, \circ\right)$, is called the algebra of maps of $Tr_{\mathcal{A}}$, on himself. 

A map $\alpha:Tr_{\mathcal{A}} \longrightarrow Tr_{\mathcal{A}}$, is an endomorphism of the group $\left(Tr_{\mathcal{A}}, \circ\right)$, on himself (see \cite{3}, \cite{15}), namely such that,
\begin{equation} \label{17}
\forall \sigma_{1}, \sigma_{2} \in Tr_{\mathcal{A}}, \alpha (\sigma_{1} \circ \sigma_{2})=\alpha (\sigma_{1}) \circ \alpha (\sigma_{2}).
\end{equation}

\begin{lemma} \label{lem.1} 
The addition of, each two endomorphisms of $Tr_{\mathcal{A}}$ on himself, is a endomorphisms of $Tr_{\mathcal{A}}$ on himself.
\end{lemma}
\proof Se Lemma Lemma 1, pp.4 to \cite{8}  \qed

\begin{lemma} \label{lem.2}
The composition of any two endomorphisms of $Tr_{\mathcal{A}}$ on himself, is an endomorphisms of $Tr_{\mathcal{A}}$ on himself.
\end{lemma}
\proof Se Lemma Lemma 2, pp.5 to \cite{8}  \qed

From Lemmas ~\ref{lem.1} \& ~\ref{lem.2}, based on the understanding of a substructure of an algebraic structure (see \cite{Lang2002}, \cite{Rotman2010}), we get this too,
\begin{theorem} \label{thm.1} \cite{8}
The set of endomorphisms of $Tr_{\mathcal{A}}$ on himself, regarding actions 'addition $+$' and 'composition $\circ$' in ti, is a substructure of algebra $\left({(Tr_{\mathcal{A}})}^{Tr_{\mathcal{A}}},+, \circ\right)$ of maps, of $Tr_{\mathcal{A}}$ on oneself.
\end{theorem}

We call it, the '\textit{endomorphisms-algebra}' of the $Tr_{\mathcal{A}}$ group on ourselves and mark with $End_{Tr_{\mathcal{A}}}$.

\begin{definition} \cite{8}
The 'tracer-preserving' endomorphism of the group $(Tr_{\mathcal{A}}, \circ)$ above itself, it is called an endomorphism $\alpha \in End_{Tr_{\mathcal{A}}}$ his, such that
\begin{equation} \label{18}
\forall \sigma \in Tr_{\mathcal{A}}, \pi_{\alpha(\sigma)}=\pi_{\sigma},
\end{equation}
otherwise, any trace according to $\alpha(\sigma)$ is a trace according to $\sigma$.
\end{definition}
The set of 'tracer-preserving' endomorphisms, of the group $(Tr_{\mathcal{A}}, \circ)$ above itself, we will marked with, $End_{Tr_{\mathcal{A}}}^{TP}$

The map $0_{Tr_{\mathcal{A}}}: Tr_{\mathcal{A}} \longrightarrow Tr_{\mathcal{A}}$, defined by

\begin{equation} \label{19}
0_{Tr_{\mathcal{A}}}(\sigma)=id_{\mathcal{P}}, \forall \sigma \in Tr_{\mathcal{A}}.
\end{equation}
is an endomorphism of the translation group $Tr_{\mathcal{A}}$ on himself \cite{8}. 

\begin{proposition} \cite{8}
The \textbf{zero endomorphism} $0_{Tr_{\mathcal{A}}}$ of $Tr_{\mathcal{A}}$ on himself, is a trace-preserving endomorphism of $Tr_{\mathcal{A}}$ on himself. 
\end{proposition}

The identical map $1_{Tr_{\mathcal{A}}}: Tr_{\mathcal{A}} \longrightarrow Tr_{\mathcal{A}}$, defined by
\begin{equation} \label{20}
1_{Tr_{\mathcal{A}}}(\sigma)=\sigma, \forall \sigma \in Tr_{\mathcal{A}},
\end{equation}

is also, an endomorphism of the translation group $Tr_{\mathcal{A}}$ on himself, see \cite{8}.

\begin{proposition} \cite{8}
The \textbf{unitary endomorphism} '$1_{Tr_{\mathcal{A}}}$' of $Tr_{\mathcal{A}}$ on himself, is a trace-preserving endomorphism of $Tr_{\mathcal{A}}$ on himself.
\end{proposition}

\begin{theorem} \label{thm.2} 
If $\alpha$ and $\beta$ are two trace-preserving endomorphisms of $Tr_{\mathcal{A}}$ on himself, then their sum $\alpha + \beta$ is an trace-preserving endomorphism.
\end{theorem}

\proof See Theorem 5, pp.6 to \cite{8}.  \qed

\begin{theorem} \label{thm.3}
If $\alpha$ and $\beta$ are two trace-preserving endomorphisms of $Tr_{\mathcal{A}}$ on himself, then their composition $\alpha \circ \beta$ is an trace-preserving endomorphism.
\end{theorem}

\proof See Theorem 6, pp.6-7 to \cite{8}.  \qed

We note now, with $\varphi:Tr_{\mathcal{A}}\longrightarrow Tr_{\mathcal{A}}$, defined by
\begin{equation} \label{21}
\forall \sigma \in Tr_{\mathcal{A}}, \varphi(\sigma)={\sigma}^{-1}
\end{equation}
is an endomorphism of the commutative group of translations $Tr_{\mathcal{A}}$ on himself, because

\begin{align*}
\forall \sigma_{1}, \sigma_{2} \in Tr_{\mathcal{A}}, \varphi(\sigma_{1} \circ \sigma_{2}) &= {(\sigma_{1} \circ \sigma_{2})}^{-1} \ (by, \ eq. \ref{21}) \\
&= {(\sigma_{2} )}^{-1} \circ {(\sigma_{1})}^{-1}  \\
&= {(\sigma_{1} )}^{-1} \circ {(\sigma_{2})}^{-1}  \  (by, \ref{14})\\
&= \varphi(\sigma_{1}) \circ \varphi(\sigma_{2}) 
\end{align*}

Also, we have that,
\[
\forall \sigma \in Tr_{\mathcal{A}}, \pi_{\sigma}=\pi_{{\sigma}^{-1}} \Rightarrow \pi_{\varphi(\sigma)}=\pi_{{\sigma}^{-1}}= \pi_{\sigma}.
\]

from this, we have prove this

\begin{proposition}
The endomorphisms $\varphi:Tr_{\mathcal{A}}\longrightarrow Tr_{\mathcal{A}}$, defined by $\varphi(\sigma)={\sigma}^{-1},\forall \sigma \in Tr_{\mathcal{A}},$ is an trace-preserving endomorphisms of $Tr_{\mathcal{A}}$ on himself.
\end{proposition}

\begin{proposition} \cite{8}
For a trace-preserving endomorphism $\alpha \in End_{Tr_{\mathcal{A}}}^{TP}$, the map $-\alpha: Tr_{\mathcal{A}}\longrightarrow Tr_{\mathcal{A}}$, defined by
\begin{equation} \label{21'}
\forall \sigma \in Tr_{\mathcal{A}}, (-\alpha)(\sigma)=(\alpha(\sigma))^{-1} \Rightarrow (-\alpha) \in End_{Tr_{\mathcal{A}}}^{TP}
\end{equation}
so,  $-\alpha$ is an trace-preserving endomorphism of $Tr_{\mathcal{A}}$ on himself.
\end{proposition}

\proof See Proposition 5, pp.7 in \cite{8}.  \qed

The endomorphism $'-\alpha'$ we call it the \textbf{additive inverse} endomorphism of endomorphism $\alpha$.

Let it be $\delta$ an dilation. According to theorem \ref{thm.normal}, have that,

\[ \forall \sigma \in Tr_{\mathcal{A}}, {\delta}^{-1} \circ \sigma \circ \delta \in Tr_{\mathcal{A}}. \]

For an dilation $\delta \in Dil_{\mathcal{A}}$, the map $\alpha_{\delta}:Tr_{\mathcal{A}}\longrightarrow Tr_{\mathcal{A}}$, defined by

\begin{equation} \label{22}
\alpha_{\delta}(\sigma)={\delta}^{-1} \circ \sigma \circ \delta, \forall \sigma \in Tr_{\mathcal{A}},
\end{equation}
is an endomorphisms of translation group $Tr_{\mathcal{A}}$ on himself, because

\begin{align*}
\forall \sigma_{1}, \sigma_{2} \in Tr_{\mathcal{A}}, \alpha_{\delta}(\sigma_{1} \circ \sigma_{2}) &= \delta^{-1} \circ (\sigma_{1} \circ \sigma_{2}) \circ \delta  \\
&= \delta^{-1} \circ \sigma_{1} \circ (\delta \circ {\delta}^{-1}) \circ \sigma_{2} \circ \delta \\
&=({\delta}^{-1} \circ \sigma_{1} \circ \delta) \circ ({\delta}^{-1} \circ \sigma_{2} \circ \delta) \\
&= \alpha_{\delta}(\sigma_{1}) \circ \alpha_{\delta}(\sigma_{2}) 
\end{align*}

From corollary \ref{cor.1} of Theorem \ref{thm.normal} we have
\[
\forall \sigma \in Tr_{\mathcal{A}}, \pi_{\sigma}=\pi_{{\delta}^{-1} \circ \sigma \circ \delta} \Rightarrow \pi_{\alpha_{\delta}}=\pi_{\sigma}.
\]

from this, we have prove this

\begin{proposition} \label{alpha.delta}
For an dilation $\delta \in Dil_{\mathcal{A}}$, the endomorphism $\alpha_{\delta}:Tr_{\mathcal{A}}\longrightarrow Tr_{\mathcal{A}}$, defined by $\alpha_{\delta}(\sigma)={\delta}^{-1} \circ \sigma \circ \delta, \forall \sigma \in Tr_{\mathcal{A}},$ is an trace-preserving endomorphisms of $Tr_{\mathcal{A}}$ on himself.
\end{proposition}

\section{The Skew-Field of Trace-Preserving Endomorphisms}

Consider now the set 
\begin{equation} \label{23}
End_{Tr_{\mathcal{A}}}^{TP}=\left\{\alpha \in End_{Tr_{\mathcal{A}}} \ | \ \alpha - \text{is a trace-preserving endomrphism}\right\}
\end{equation}
of trace-preserving endomorphisms of $Tr_{\mathcal{A}}$ in itself. According to Theorem \ref{thm.2} and Theorem \ref{thm.3}, $(End_{\mathcal{A}}^{TP}, +, \circ)$ is a substructure of the algebra  $(End_{\mathcal{A}}, +, \circ)$, of endomorphisms related to addition and composition actions, therefore it is itself an algebra.

\begin{theorem} \label{additive.group} \cite{8}
The Grupoid $\left(End_{Tr_{\mathcal{A}}}^{TP}, +\right)$, is commutative (Abelian) Group.
\end{theorem}

\proof  See Theorem 7, pp.7-8 to \cite{8}.  \qed

\begin{theorem} \label{unitary.ring} \cite{8}
The algebra $\left(End_{Tr_{\mathcal{A}}}^{TP}, + , \circ\right)$, is a Associative Unitary Ring.
\end{theorem}

\proof  See Theorem 8, pp.8-9 to \cite{8}.  \qed

For attaining, the algebraic structure we aim for, namely the skew-field of trace-preserving endomorphism, we need that semi-group $(End_{Tr_{\mathcal{A}}} -\left\{0_{Tr_{\mathcal{A}}}\right\}, \circ)$ to become a Group, of course firstly the set $End_{Tr_{\mathcal{A}}} -\left\{0_{Tr_{\mathcal{A}}}\right\}$ should be closed in relation to the composting action (This is proved in the theorem \ref{skew.field}). 
We are first giving the following Lemma, which will help us to find the inverse of a trace-preserving endomorphisms related to the composition action (analogy of the multiplication action) for the algebraic structure we are building. So any trace-preserving endomorphism different from zero endomorphism has its inverse. 

\begin{lemma} \label{lem.invers}
For all $\alpha \in End_{Tr_{\mathcal{A}}}-\left\{0_{Tr_{\mathcal{A}}}\right\}$, exists, only an dilation $\delta \in Dil_{\mathcal{A}}$, of the affine plane that even has a fixed-point regarding $\delta$, such that
\begin{equation} \label{24}
\forall \sigma \in Tr_{\mathcal{A}}, \ \alpha(\sigma)=\delta \circ \sigma \circ {\delta}^{-1}.
\end{equation}
\end{lemma}

\proof For every two points $R,Q$ of the affine plane $\mathcal{A}$ always, exist a translation that brings point $R$ to point $Q$. If $Q \neq R,$ according to axiom A.1 of affine plane, has a single line $\ell^{RQ}$ in this plane. If $Q=R,$ such is the identical translation $id_{\mathcal{A}}$; even in this case, a line passing through the point Q will mark $\ell^{QQ}$. So exist and is only a translation of the affine plans, Which we mark with $\sigma_{RQ}$, such that,
\begin{equation} \label{25}
\sigma_{RQ}(R)=Q
\end{equation}
We show at the beginning that, if there exists the dilation $\delta \in Dil_{\mathcal{A}}$, satisfying eq. \ref{24}, with fixed-point for it, the point $P$, then he is alone.
Let $Q$ be a point of affine plane $\mathcal{A}$ and $\sigma_{PQ}$ the translation that $\sigma_{PQ}(P)=Q$. Then,
\[  \alpha(\sigma_{PQ})=\delta \circ \sigma_{PQ} \circ {\delta}^{-1}. \]
Meanwhile,
\[ \delta(P)=P \Longleftrightarrow \delta^{-1} (P)=P, \]
so point $P$ is a fixed-point even by $\delta^{-1}$. Under these conditions the image of point $P$ by translation $\alpha(\sigma_{PQ})$, is
\begin{align*}
 [\alpha(\sigma_{PQ})] (P) &= [\delta \circ \sigma_{PQ} \circ {\delta}^{-1}] (P)  \\
&= [\delta \circ \sigma_{PQ}]({\delta}^{-1} (P))  \\
&= [\delta \circ \sigma_{PQ}](P)   \\
&= \delta (\sigma_{PQ}(P) ) \\
&= \delta (Q)  \ \ (by, \ eq. \ref{25}).
\end{align*}

So, for fixed point $P$ related with dilation  $\delta$, and $\forall Q \in \mathcal{P}$ in affine plane $\mathcal{A}$, have
\begin{equation} \label{26}
\delta (Q)=[\alpha(\sigma_{PQ})](P).
\end{equation}

So if exist an dilation  satisfying eq. \ref{24}, being $P$ fixed point on it, then his image of a point $Q$ of the affine plane, is given by eq. \ref{26} as follows: constructed the translation $\sigma_{PQ}$, from it, constructed the translation $\alpha(\sigma_{PQ})$, and the latter is applied to the point $P$. By this we showed that dilation $\delta$ is defined in a single way.

Let it be now $\alpha \in End_{Tr_{\mathcal{A}}}-\left\{0_{Tr_{\mathcal{A}}}\right\}$ and define a map $\delta$ of the affine plane according to eq. \ref{26}. Let be $Q, R$ two different points of the affine plane $\mathcal{A}$, for translations $\sigma_{QR}$ and $\sigma_{PQ}$, from eq.\ref{25} we have,
\[ [\sigma_{QR} \circ \sigma_{PQ}](P)=\sigma_{QR}(Q)=R. \]

This indicates that the translation $[\sigma_{QR} \circ \sigma_{PQ}]$, leads point $P$ to point $R$, as well as $\sigma_{PR}$. From here, by Corollary 3.4. in \cite{5} have that, $[\sigma_{QR} \circ \sigma_{PQ}]=\sigma_{PR}$, then

\begin{align*}
 [\alpha(\sigma_{PR})] (P) &= [\alpha(\sigma_{QR} \circ \sigma_{PQ})](P)  \\
&= [\alpha(\sigma_{QR}) \circ \alpha(\sigma_{PQ})](P) \ \ (by, \ eq. \ref{17}) \\
&= [\alpha(\sigma_{QR})]  \left(\alpha(\sigma_{PQ})(P) \right)  .
\end{align*}

So,
\[
[\alpha(\sigma_{QR})]  \left(\circ \alpha(\sigma_{PQ})(P) \right)=[\alpha(\sigma_{PR})] (P),
\]
from this we have that,
\begin{equation} \label{27}
[\alpha(\sigma_{QR})]  \left(\delta(Q) \right)=\delta(R),
\end{equation}
because, according to eq. \ref{26}, along with $\left(\alpha(\sigma_{PQ})(P) \right)=\delta(Q)$, we have and $[\alpha(\sigma_{PR})](P)=\delta(R)$.

Consider now the line $\ell_{QR}^{\delta(Q)}$, that passes by the point $\delta(Q)$ and is parallel with line $QR=Q\sigma_{QR}(Q)$. So it is a trace of the point $\delta(Q)$, according to translation $\sigma_{QR}$. But  
\[ \alpha \neq 0_{Tr_{\mathcal{A}}}\Rightarrow \alpha(\sigma_{QR}) \neq id_{\mathcal{A}}, \]
consequently has a certain direction, that defines a trace, in particular the trace of the point $ \delta(Q)$ according to it. Meanwhile, it is a trace-preserving endomorphism, therefore this trace is the line $\ell_{QR}^{\delta(Q)}$. So we have that $\delta(R) \in \ell_{QR}^{\delta(Q)}$. 
In this way we have prove that for every two different points $Q, R$  of the affine plane we have

\begin{equation} \label{28}
\ell^{QR} \parallel \ell^{\delta(Q)\delta(R)}.
\end{equation}
In addition, for map $\delta$ we also have:
\begin{align*}
\delta(P) &= [\alpha(\sigma_{PP})](P)  \ \ (by, \ eq. \ref{26}) \\
&= [\alpha(id_{\mathcal{A}})](P)  \\
&= ( id_{\mathcal{A}} (P) )\\
&= P,
\end{align*}
that brings,
\begin{equation} \label{29}
\delta(P)= P.
\end{equation}

From a theorem in \cite{11}, a map of affine planes, which for every two different points $P, Q$ satisfies \ref{28}, or it is a degenerate map, in the sense that each point of the plan has the same image as a given point, or is a bijection. If the map $\delta$ is not a bijection, then it is such that every point $Q$ of the affine plane has the same image with a given point, namely all points of affine plane have the same image. Hence,
\[
\delta(Q)=\delta(P)=P\Rightarrow P=[\alpha(\sigma_{PQ})](P) \Rightarrow \alpha(\sigma_{PQ})=id_{\mathcal{A}}, \forall Q \in \mathcal{P}.
\]
It is clear that any translation $\sigma$ can be represented in the form $\sigma_{PQ}$. So the last statement is written in the form,
\[ \alpha(\sigma)=id_{\mathcal{A}}, \ \forall \sigma \in Tr_{\mathcal{A}}. \]

It shows, according to \ref{19}, that $\alpha=0_{\mathcal{A}}$ which is a contradiction, since $\alpha \in End_{Tr_{\mathcal{A}}}^{TP}-\left\{0_{\mathcal{A}}\right\}$. It remains that the map $\delta$ is bijection, therefore
\[
Q \neq R \Rightarrow \delta(Q) \neq \delta(R).
\]
Then (\ref{28}) shows that $\delta$ is an dilation, while \ref{29} indicates that, the point $P$ is a fixed point by $\delta$.
Finally we show that this dilation satisfies \ref{24}. For this, according to \ref{26} and \ref{28}, we have,
\[
\delta(Q)=[\alpha(\sigma_{PQ})](\delta(P))=[\alpha(\sigma_{PQ}) \circ \delta](P).
\]
From here we present the point $Q$ in the form
\[
Q=\delta^{-1}\left([\alpha(\sigma_{PQ}) \circ \delta](P)\right),
\]
from where, have
\[
Q=\left(\delta^{-1} \circ \alpha(\sigma_{PQ}) \circ \delta\right) (P),
\]

which indicates that the translation  $\delta^{-1} \circ [\alpha(\sigma_{PQ})] \circ \delta$, leads point $P$ to point $Q$. Therefore, it is the translation $\sigma_{PQ}$. Hence,

\[
\delta^{-1} \circ [\alpha(\sigma_{PQ})] \circ \delta=\sigma_{PQ} \Leftrightarrow \alpha(\sigma_{PQ})=\delta \circ \sigma_{PQ} \circ \delta^{-1}.
\]

We highlighted above that any translation $\sigma$, is presented as $\sigma_{PQ}$. So we get the assertion,

\[
\forall \sigma \in Tr_{\mathcal{A}}, \  \alpha(\sigma)=\delta \circ \sigma \circ \delta^{-1}.
\]

\qed

\begin{theorem} \label{inverse.thm}
For all $\alpha \in End_{Tr_{\mathcal{A}}}^{TP}-\left\{0_{\mathcal{A}}\right\}$, exist its inverse $\alpha^{-1} \in End_{Tr_{\mathcal{A}}}^{TP}-\left\{0_{\mathcal{A}}\right\}$.
\end{theorem}

\proof For $\alpha \in End_{Tr_{\mathcal{A}}}^{TP}-\left\{0_{\mathcal{A}}\right\}$, according to Lemma \ref{lem.invers} find the dilation $\delta$, such that it applies \ref{24}.

We construct now for this, dilation $\delta$ the map $\alpha_{\delta}$, determined by
\begin{equation} \label{*}
\alpha_{\delta}(\sigma)=\delta^{-1} \circ \sigma \circ \delta, \ \forall \sigma \in Tr_{\mathcal{A}}.
\end{equation}

According to Proposition \ref{alpha.delta}, the map $\alpha_{\delta}$ is an trace-preserving endomorphism even different from $0_{\mathcal{A}}$, therefore also $\alpha_{\delta} \in End_{Tr_{\mathcal{A}}}^{TP}-\left\{0_{\mathcal{A}}\right\}$. Then we have

\begin{align*}
\forall \sigma \in Tr_{\mathcal{A}},    [\alpha \circ \alpha_{\delta}] (\sigma) &= \alpha [\alpha_{\delta} (\sigma)]  \ (by, \ eq.\ref{16}) \\
&= \alpha [\delta^{-1} \circ \sigma \circ \delta] \ \ (by, \ eq. \ref{*}) \\
&= \delta \circ [\delta^{-1} \circ \sigma \circ \delta] \circ \delta^{-1} \  (by, \ eq. \ref{24}) \\
&= (\delta \circ \delta^{-1}) \circ \sigma \circ (\delta \circ \delta^{-1}) \\
&= \sigma \\
&= 1_{Tr_{\mathcal{A}}}.
\end{align*}

As well, we have,

\begin{align*}
\forall \sigma \in Tr_{\mathcal{A}},   [\alpha_{\delta} \circ \alpha] (\sigma) &= \alpha_{\delta} [\alpha (\sigma)] \ \ (by, \ eq. \ref{16})  \\
&= \alpha_{\delta} [\delta \circ \sigma \circ \delta^{-1}] \ \ (by, \ eq. \ref{24}) \\
&= \delta^{-1} \circ [\delta \circ \sigma \circ \delta^{-1}] \circ \delta \ \ (by, \ eq. \ref{*}) \\
&= (\delta^{-1} \circ \delta) \circ \sigma \circ (\delta^{-1} \circ \delta) \\
&= \sigma \\
&= 1_{Tr_{\mathcal{A}}}.
\end{align*}
So, we have
\[
\forall \alpha \in End_{Tr_{\mathcal{A}}}^{TP}-\left\{0_{\mathcal{A}}\right\}, \ \alpha \circ \alpha_{\delta}=\alpha_{\delta} \circ \alpha=1_{Tr_{\mathcal{A}}}
\]
So exist the inverse element, $\alpha^{-1} \in End_{Tr_{\mathcal{A}}}^{TP}-\left\{0_{\mathcal{A}}\right\}$, of the map $\alpha \in End_{Tr_{\mathcal{A}}}^{TP}-\left\{0_{\mathcal{A}}\right\}$, which is endomorphism $\alpha_{\delta}$ determined by \ref{*}.
\qed

\begin{theorem} \label{skew.field}
The Algebra $\left(End_{Tr_{\mathcal{A}}}^{TP}, +, \circ\right)$ is a Skew-Field.
\end{theorem}

\proof From the Theorem \ref{unitary.ring} (see Theorem 8, pp. 8-9, in \cite{8}) we have that, the algebra $(End_{Tr_{\mathcal{A}}}^{TP}, + , \circ)$, is a associative unitary Ring.

Also because $1_{Tr_{\mathcal{A}}} \neq 0_{Tr_{\mathcal{A}}}$, so in $End_{Tr_{\mathcal{A}}}^{TP}$ has at least one nonzero element, according to the definition of a skew-field \cite{Lang2002}, \cite{Rotman2010}, \cite{RobertWisbauer.RingTheory1991}, required to prove the following:

1. $End_{Tr_{\mathcal{A}}}^{TP}-\left\{0_{Tr_{\mathcal{A}}}\right\}$, is a closed subset of $End_{Tr_{\mathcal{A}}}^{TP}$ about the composition. Really, if $\alpha_{1} \neq 0_{Tr_{\mathcal{A}}}$ and $\alpha_{2} \neq 0_{Tr_{\mathcal{A}}}$, then, and $\alpha_{1} \circ \alpha_{2} \neq 0_{Tr_{\mathcal{A}}}$. 
Suppose that $\alpha_{1} \circ \alpha_{2} = 0_{Tr_{\mathcal{A}}}$. Given that $\alpha_{1} \neq 0_{Tr_{\mathcal{A}}}$, then exist the inverse $\alpha_{1}^{-1} \in End_{Tr_{\mathcal{A}}}^{TP}-\left\{0_{Tr_{\mathcal{A}}}\right\}$, hence 

\begin{align*}
\alpha_{1} \circ \alpha_{2} = 0_{Tr_{\mathcal{A}}}\Rightarrow  0_{Tr_{\mathcal{A}}} &= \alpha_{1}^{-1} \circ (\alpha_{1} \circ \alpha_{2})   \\
&= (\alpha_{1}^{-1} \circ \alpha_{1}) \circ \alpha_{2} \\
&= \alpha_{2}.
\end{align*}

This is contradict with the condition that $\alpha_{2} \neq 0_{Tr_{\mathcal{A}}}$. Same is to suppose that $\alpha_{2} \neq 0_{Tr_{\mathcal{A}}}$, would result that $\alpha_{1}=0_{Tr_{\mathcal{A}}}$, that is also a contradict with the condition $\alpha_{1} \neq 0_{Tr_{\mathcal{A}}}$.
So,
\[
\forall \alpha_{1}, \alpha_{2} \neq 0_{Tr_{\mathcal{A}}}, \ \alpha_{1} \circ \alpha_{2} \neq 0_{Tr_{\mathcal{A}}}.
\]
Hence, our associative unitary ring has not divisor of zero.

2. The Grupoid $\left(End_{Tr_{\mathcal{A}}}^{TP}-\left\{0_{Tr_{\mathcal{A}}}\right\}, \circ\right)$ is Group.
Composition '$\circ$', being associative in $End_{Tr_{\mathcal{A}}}^{TP}$, then he is associative in $End_{Tr_{\mathcal{A}}}^{TP}-\left\{0_{Tr_{\mathcal{A}}}\right\}$. In addition, given that $1_{Tr_{\mathcal{A}}} \in End_{Tr_{\mathcal{A}}}^{TP}$ then he is the unitary element of $\left(End_{Tr_{\mathcal{A}}}^{TP}-\left\{0_{Tr_{\mathcal{A}}}\right\}, \circ\right)$.
According to Theorem \ref{inverse.thm} exist inverse element of $\alpha^{-1}$, for all $\alpha \in End_{Tr_{\mathcal{A}}}^{TP}-\left\{0_{Tr_{\mathcal{A}}}\right\}$, hence, $\left(End_{Tr_{\mathcal{A}}}^{TP}-\left\{0_{Tr_{\mathcal{A}}}\right\}, \circ\right)$ is Group.

\qed


\bibliographystyle{amsplain}
\bibliography{NSrefs}

\end{document}